\documentclass[10pt]{amsart}
\usepackage{amsmath, amssymb, amsthm, a4wide}
\usepackage{geometry}
\usepackage{fancyhdr}
\usepackage{array}
\usepackage{enumerate}
\newtheorem{theorem}{Theorem}[section]
\newtheorem{lemma}[theorem]{Lemma}
\newtheorem{definition}[theorem]{Definition}

\newtheorem{remark}[theorem]{Remark}
\newtheorem{assumption}[theorem]{Assumption}
\numberwithin{equation}{section}

\newcounter{newlist} 
\newenvironment{mylist}[1][]{
\begin{list}{\textnormal{(\arabic {newlist})}} 
    { 
    \usecounter{newlist} 
     \setlength{\labelsep}{0.5em}
     \setlength{\leftmargin}{1.7em}
     \setlength{\rightmargin}{0cm} 
     \setlength{\parsep}{0pt}
     \setlength{\itemsep}{3pt}
     \setlength{\itemindent}{0em} 
     \setlength{\listparindent}{2em} 
    }}
{\end{list}}
\newcounter{nnnewlist} 
\newenvironment{mmmylist}[1][]{
\begin{list}{\textnormal{(\arabic {nnnewlist})}} 
    { 
    \usecounter{nnnewlist} 
     \setlength{\labelsep}{0.5em}
     \setlength{\leftmargin}{1.8em} 
     \setlength{\rightmargin}{0cm} 
     \setlength{\parsep}{0pt}
     \setlength{\itemsep}{3pt} 
     \setlength{\itemindent}{-1.3em}
     \setlength{\listparindent}{2em} 
    }}
{\end{list}}
\newcounter{nelist} 

\newcounter{nmelist} 
\newenvironment{mnlist}[1][]{
\begin{list}{\textnormal{(H\arabic {nmelist})}} 
    { 
    \usecounter{nmelist} 
     \setlength{\labelsep}{0.5em} 
     \setlength{\leftmargin}{0.5em} 
     \setlength{\rightmargin}{0cm} 
     \setlength{\parsep}{0pt} 
     \setlength{\itemsep}{3pt} 
     \setlength{\itemindent}{0em} 
     \setlength{\listparindent}{2em} 
    }}
{\end{list}}
\newcounter{nnmelist} 
\newenvironment{mnnlist}[1][]{
\begin{list}{\textnormal{(H\arabic {nnmelist}')}} 
    { 
    \usecounter{nnmelist} 
     \setlength{\labelsep}{0.5em} 
     \setlength{\leftmargin}{0.5em} 
     \setlength{\rightmargin}{0cm} 
     \setlength{\parsep}{0pt}
     \setlength{\itemsep}{3pt}
     \setlength{\itemindent}{0em} 
     \setlength{\listparindent}{2em} 
    }}
{\end{list}}
\newcounter{nnmnelist} 
\newenvironment{nmnnlist}[1][]{
\begin{list}{\textnormal{(A\arabic {nnmnelist})}} 
    { 
    \usecounter{nnmnelist} 
     \setlength{\labelsep}{0.5em}
     \setlength{\leftmargin}{0.5em} 
     \setlength{\rightmargin}{0cm} 
     \setlength{\parsep}{0pt}
     \setlength{\itemsep}{3pt} 
     \setlength{\itemindent}{0em}
     \setlength{\listparindent}{2em} 
    }}
{\end{list}}
\newcounter{selist} 
\newenvironment{slist}[1][]{
\begin{list}{\textnormal{(H\arabic {selist}")}} 
    { 
    \usecounter{selist} 
     \setlength{\labelsep}{0.5em} 
     \setlength{\leftmargin}{3em} 
     \setlength{\rightmargin}{0cm} 
     \setlength{\parsep}{0pt} 
     \setlength{\itemsep}{3pt} 
     \setlength{\itemindent}{0em} 
     \setlength{\listparindent}{2em} 
    }}
{\end{list}}

\begin{document}

\title[GSDEs with Integral-Lipschitz Coefficients]{
On the Existence and Uniqueness of Solutions\\to
Stochastic Differential Equations\\Driven by $G$-Brownian
Motion\\with Integral-Lipschitz Coefficients}
\author{Xuepeng BAI}
\author{Yiqing LIN}
\address{School of Mathematics\newline\indent Shandong University\newline\indent 250100 Jinan, China\newline\newline\indent Institut de Recherche Math\'{e}matique de Rennes\newline\indent Universit\'{e} de Rennes 1\newline\indent 35042 Rennes Cedex, France\newline}
\email{yiqing.lin@univ-rennes1.fr (Corresponding author: Yiqing LIN)}
\thanks{The research is partially supported by the Major Program in Key Research Institute of Humanities and Social Sciences sponsored by Ministry of Education of China under grant No. 2009JJD790049 and the Post-graduate Study Abroad Program sponsored by China Scholarship Council.}


\date{November 1st, 2012}


\begin{abstract}In this paper, we study
the existence and uniqueness of solutions to stochastic differential
equations driven by $G$-Brownian motion (GSDEs) with integral-Lipschitz
conditions on their coefficients.
\end{abstract}

\maketitle

\section{Introduction}
\noindent Motivated by uncertainty problems, risk
measures and super-hedging in finance, Peng \cite{P2, P3} 
introduced a framework of $G$-expectation, in which a new type of Brownian motion was constructed and the related stochastic calculus has been established. As a counterpart in the classical framework, stochastic differential equations driven by $G$-Brownian motion (GSDEs) have been studied by Gao \cite{G} and Peng \cite{P3}. In these works, the solvability of GSDEs under Lipschitz conditions has been obtained by contraction mapping theorem. \\[6pt] 
\noindent Typically, a GSDE is of the following form:
\begin{equation}\label{eq 2}
X(t)=x+\int_0^tb(s,X(s))ds+\int_0^th(s,X(s))d\langle
B,B\rangle_s+\int_0^tg(s,X(s))dB_s,\ 0\leq t\leq T,
\end{equation}
where $x\in\mathbb{R}^n$ is the initial value, $B$ is the $G$-Brownian motion and $\langle B,B\rangle$ is the quadratic variation process of $B$.\\[6pt]
\noindent In this paper, we study the solvability of the GSDE (\ref{eq 2}) under a so-called integral-Lipschitz condition:
\begin{equation}\label{eq 3}
|b(t,x_1)-b(t,x_2)|^2+|h(t,x_1)-h(t,x_2)|^2+|g(t,x_1)-g(t,x_2)|^2\leq\rho(|x_1-x_2|^2),
\end{equation}
where  $\rho:(0,+\infty)\rightarrow(0,+\infty)$ is a continuous
increasing and concave function that vanishes at $0+$ and satisfies
$$\int_0^1\frac{dr}{\rho(r)}=+\infty.$$
A typical example of (\ref{eq 3}) is
\begin{align*}
|b(t,x_1)-b(t,x_2)|^2+|h(t,x_1)-h(t,x_2)|^2+|g(t,x_1)-g(t,x_2)|^2\leq|x_1-x_2|^2\ln\frac{1}{|x_1-x_2|}.
\end{align*}
Furthermore, we consider the GSDE (\ref{eq 2}) under a ``weaker'' condition on $b$ and $h$:
\begin{equation}\label{eq plus}
|b(t, x_1)-b(t, x_2)|+|h(t, x_1)-h(t, x_2)|\leq\rho(|x_1-x_2|),
\end{equation}
where $\rho$ satisfies the same conditions as above. A typical example of (\ref{eq plus}) is
\begin{align*}
|b(t, x_1)-b(t, x_2)| + |h(t, x_1)-h(t, x_2)|&\leq|x_1-x_2|\ln\frac{1}{|x_1-x_2|}.
\end{align*}
\noindent In the classical framework, Watanabe and Yamada \cite{WY, WY1} and Fang and Zhang \cite{FZ} have proved the pathwise uniqueness of solutions to finite-dimensional SDEs under some similar non-Lipschitz condition. In addition to that, Yamada \cite{Y} has found an explicit way to construct the solutions by successive approximation. On the other hand, Hu and Lerner \cite{HL} have worked on the SDEs in infinite dimension under the integral-Lipschitz conditions (\ref{eq 3}) and (\ref{eq plus}). They established both the pathwise uniqueness and successive approximations of the solutions. 
Corresponding to the result in Watanabe and Yamada \cite{WY}, Lin \cite{LQ} has obtain a pathwise uniqueness result for non-Lipschitz GSDEs when the coefficient $g$ is bounded.\\[6pt]
In this article, we present both the existence and uniqueness results for GSDE (\ref{eq 2}) under the integral-Lipschitz conditions (\ref{eq 3}) and (\ref{eq plus}). These results are obtained by a similar technique adopted by Hu and Lerner \cite{HL}. This paper is organized as follows: Section 2 gives the necessary preliminaries in the $G$-framework.
Section 3 proves the existence and uniqueness theorem for GSDEs with integral-Lipschitz coefficients and Section 4 studies the case for $G$-backward stochastic differential equations (GBSDEs).
\section{Preliminaries}
\noindent The main purpose of this section is to recall some preliminary results in the $G$-framework, which are necessary later in the text. The reader interested
in a more detailed description of these notions is referred to Denis et al. \cite{DHP}, Gao \cite{G} and Peng \cite{P3}.
\subsection{$G$-Brownian motion and $G$-expectation}
Adapting the approach in Peng \cite{P3}, let
$\Omega$ be a given nonempty fundamental space and $\mathcal {H}$ a linear space of real
functions defined on $\Omega$ such that
\begin{mylist}
\item $1\in \mathcal {H}$.
\item $\mathcal {H}$ is stable with respect to bounded Lipschitz
functions, i.e., for all $n\geq1$, $X_1,\ldots,X_n\in
\mathcal {H}$ and $\varphi \in C_{b, lip}(\mathbb{R}^n)$, it holds also
$\varphi(X_1,\ldots,X_n)\in\mathcal {H}$.
\end{mylist}
\begin{definition}\label{de
2.1.1}A sublinear expectation $\mathbb{E}[\cdot]$ on $\mathcal {H}$ is a
functional $\mathbb{E}[\cdot]:\mathcal {H}\rightarrow\mathbb{R}$ with the following
properties: for
each $X$, $Y\in \mathcal {H}$, we have
\begin{mmmylist}
\item \textbf{Monotonicity:} if $X\geq Y$, then $\mathbb{E}[X]\geq\mathbb{E}[Y];$
\item \textbf{Preservation of constants:}\ $\mathbb{E}[c]=c,$ for all $c\in \mathbb{R}$;
\item \textbf{Sub-additivity:}\ $\mathbb{E}[X]-\mathbb{E}[Y]\leq\mathbb{E}[X-Y]$;
\item \textbf{Positive homogeneity:}\ $\mathbb{E}[\lambda
X]=\lambda\mathbb{E}[X]$, for all $\lambda\in\mathbb{R}^+$.
\end{mmmylist}
\end{definition}
The triple $(\Omega,\mathcal {H}, \mathbb{E})$  is called a
sublinear expectation space. 
\begin{definition}\label{de 2.1.2} A random vector
$Y=(Y_1,\ldots,Y_n) \in \mathcal {H}^n$ is said to be
independent of $X\in\mathcal {H}^m$ under $\mathbb{E}[\cdot]$ if for each test function $\varphi\in C_{b, lip}(\mathbb{R}^{n+m}) $ we
have
$$\mathbb{E}[\varphi(X,Y)]=\mathbb{E}[\mathbb{E}[\varphi(x,Y)]_{x=X}].$$
\end{definition}
\begin{definition}
Let $X=(X_1,\ldots,X_n)\in\mathcal {H}^{ n}$ be a given random vector. We define 
$$\mathbb{F}_X[\varphi]:=\mathbb{E}[\varphi(X)],\
\varphi\in
C_{b,lip}(\mathbb{R}^n).$$
Then, the functional $\mathbb{F}_X[\cdot]$ is called the distribution of $X$ under $\mathbb{E}[\cdot]$.
\end{definition}
\noindent Now we begin to introduce the definition of  $G$-Brownian motion and
$G$-expectation.
\begin{definition}\label{de 2.1.4}
A $d$-dimensional random vector $X$ in a sublinear expectation space
$(\Omega,\mathcal {H},\mathbb{E})$ is called $G$-normal distributed if for each $\varphi \in C_{b,lip} (\mathbb{R}^d)$,
$$
u(t,x):=\mathbb{E}[\varphi(x+\sqrt{t}X)],\ t\in\mathbb{R}^+,\ x\in \mathbb{R}^d,
$$
is the viscosity solution to the following PDE defined on
$\mathbb{R}^+\times \mathbb{R}^d$:
$$
\left\{\begin{aligned}
&\frac{\partial u}{\partial t}-G(D^2u)=0;\\
&u|_{t=0} = \varphi,
\end{aligned}\right.
$$
where $G=G_X(A):\mathbb{S}^d \rightarrow \mathbb{R}$ is defined by
$$
G_X(A):=\frac{1}{2} \mathbb{E} [(AX,X)]
$$
and $D^2 u=(\partial^2_{x_i x_j} u)^d_{i,j=1}$.
\end{definition}
\noindent In particular, $\mathbb{E}[\varphi(X)]=u(1,0)$ defines the distribution of $X$.  By Theorem 2.1 in Chapter I of Peng \cite{P3}, 
there exists a bounded and closed subset $\Gamma$ of $\mathbb{R}^d$, such
that for each $A\in \mathbb{S}^d$, $G_X(A)$ can be represented
as
$$
G_X(A)=\frac{1}{2}
\sup_{\gamma\in\Gamma} {\rm tr}[\gamma \gamma^{\bf{Tr}} A].
$$
Defining a subset $\Sigma :=\{\gamma \gamma^{\bf{Tr}}: \gamma\in\Gamma\}$ in $\mathbb{S}^d$, the $G$-normal distribution can be denoted by $\mathcal{N}(0,\Sigma)$.\\[6pt]
Let $\Omega$ be the space of all $\mathbb{R}^d$-valued
continuous paths $(\omega_t)_{t\geq0}$ that start from $0$ and $B$ the canonical process. We assume additionally that $\Omega$ is a metric  space equipped
with the following distance:
$$
\rho(\omega^1,\omega^2):=\sum^\infty_{N=1} 2^{-N} (\max_
{0\leq t\leq N} | \omega^1_t-\omega^2_t|)  \wedge 1).
$$
For a fixed $T\geq0$, we set
$$
L_{ip}^0 (\Omega_T):=\{\varphi(B_{t_1}, \ldots, B_{t_n}):n\geq
1,\ 0\leq t_1\leq \ldots\leq t_n \leq T,\ \varphi \in C_{b,lip}
(\mathbb{R}^{d\times n})\}.
$$
\begin{definition}\label{de 2.1.5}
Let $\mathbb{E}[\cdot]: L_{ip}^0(\Omega_T) \rightarrow \mathbb{R}$ be
a sublinear expectation on $L_{ip}^0(\Omega_T)$, we call
$\mathbb{E}[\cdot]$ a $G$-expectation if the canonical process
$B$ is a $G$-Brownian motion under
$\mathbb{E}[\cdot]$, that is, for each $0\leq s\leq t\leq T$, the increment $B_t - B_s$ is
$\mathcal{N}(0, (t-s)\Sigma)$-distributed and is independent of $(B_{t_1},\ldots,B_{t_n})$, for all $n\geq 1$ and $0\leq
t_1\leq\ldots\leq t_n\leq s$.
\end{definition}
\noindent We denote by $L^p_G(\Omega_T)$ the completion of $L_{ip}^0
(\Omega_T)$ under the Banach
norm $\mathbb{E}[|\cdot|^p]^\frac{1}{p},\ 1\leq p<+\infty,$ and we still use the notation $\mathbb{E}[\cdot]$ to 
denote the extension of this sublinear expectation.
\begin{definition}
Let $\mathbb{E}[\cdot]: L_{ip}^0(\Omega_T) \rightarrow \mathbb{R}$ be
a $G$-expectation on $L_{ip}^0(\Omega_T)$, we define the related conditional expectation of $X\in L^0_{ip}(\Omega_T)$ under $L^0_{ip}(\Omega_{t_j})$, $0\leq t_1\leq\ldots\leq t_{j-1}\leq t_j\leq t_{j+1}\leq\ldots\leq t_n\leq T$:
\begin{align*}
\mathbb{E}[X|\Omega_{t_j}]&:=\mathbb{E}[\varphi(B_{t_1},
\ldots,B_{t_n}-B_{t_{n-1}})|\Omega_{t_j}]=\mathbb{E}[\psi(B_{t_1},\ldots,B_{t_j}-B_{t_{j-1}})],
\end{align*}
where $\psi(x_1,\ldots,x_j):=\mathbb{E}[\varphi(x_1,\ldots,x_j,B_{t_{j+1}}-B_{t_j},\ldots, B_{t_n}-B_{t_{n-1}})].$ Moreover, the mapping
$\mathbb{E}[\cdot|\Omega_{t_j}]
:L_{ip}^0(\Omega_T)\rightarrow L_{ip}^0(\Omega_{t_j})$
can be continuously extended to $\mathbb{E}[\cdot|\Omega_{t_j}]:L_G^1(\Omega_T)\rightarrow L_G^1(\Omega_{t_j})$.
\end{definition}
\subsection{G-capacity}
Derived in Denis et al. \cite{DHP}, $G$-expectation can be formulated as an upper expectation of a weakly compact family of probability measures. This family is related to the set $\Gamma$ mentioned in the last subsection, which is a bounded and closed subset of $\mathbb{R}^d$ that characterizes the $G$-function $G(\cdot)$.\\[6pt]
Let $\mathbb{P}_0$ be the Wiener measure on $\Omega$, $\mathcal{F}$ the filtration
generated by the canonical process $B$ and 
$\mathcal {A}_{[0,+\infty)}^\Gamma$ the collection of all
$\Gamma-$valued progressively measurable processes. 
For each $\theta\in\mathcal {A}_{[0,+\infty)}^\Gamma$, let
$\mathbb{P}_\theta$ be the probability measure introduced by the following strong formulation:
$$
\mathbb{P}_\theta:=\mathbb{P}_0\circ(X_\theta)^{-1}, 
$$ 
where $X_\theta:=(\int_0^t\theta_sdB_s)_{t\geq 0}$, $\mathbb{P}_0$-a.s.. We set $\mathcal {P}:=\{\mathbb{P}_\theta: \theta\in\mathcal {A}_{[0,+\infty)}^\Gamma\}$ and denote by $\mathcal{P}_G$ the closure of $\mathcal{P}$ under the topology of weak convergence. \\[6pt]
Consider a capacity formulated by upper probability:
$$
\bar{C}(A):=\sup_{\mathbb{P}\in\mathcal{P}_G}\mathbb{P}(A),\ A\in\mathcal {B}(\Omega).
$$
By Proposition 50 in Denis et al. \cite{DHP}, 
$\mathcal {P}_G$ is weakly compact and thus,
$\bar{C}(\cdot)$ is a Choquet capacity. 
Then, we have the following notion of ``quasi-surely"(q.s.).
\begin{definition}\label{de 2.1.6}
A set $A\in\mathcal{B}(\Omega)$ is called polar if $\bar{C}(A)=0.$ A property
is said to hold quasi-surely if it holds outside a polar
set.
\end{definition}
\noindent On the other hand, we set for each $X\in L^0(\Omega_T)$,
\begin{equation}\label{grep}
\bar{\mathbb{E}}[X]:=\sup_{\mathbb{P}\in\mathcal{P}_G}E^{\mathbb{P}}[X].
\end{equation}
In (\ref{grep}), $E^\mathbb{P}[X]$ exists under each $\mathbb{P}\in\mathcal{P}_G$, so $\bar{\mathbb{E}}[X]$ is well defined. By Theorem 52 in Denis et al. \cite{DHP}, this upper expectation $\bar{\mathbb{E}}[\cdot]$ is consistent with $G$-expectation $\mathbb{E}[\cdot]$ on $L_G^1(\Omega_T)$, i.e.,  
$$
\bar{\mathbb{E}}[X]=\mathbb{E}[X],\ {\rm{for\ all}}\ X\in L^1_G(\Omega_T).
$$
Thus, from now on, we do not distinguish these two notations $\mathbb{E}[\cdot]$ and $\bar{\mathbb{E}}[\cdot]$. \\[6pt]
\noindent By the definitions of $\bar{\mathbb{E}}[\cdot]$ and $\bar{C}(\cdot)$, we can easily deduce the following Markov inequality and the upwards monotone convergence theorem in the $G$-framework:
\begin{lemma}\label{markov}
Let $X\in L^0(\Omega_T)$ and for some $p>0$, $\bar{\mathbb{E}}[|X^p|]<+\infty$. Then, for each $M>0$, $$
\bar{C}(\{|X|>M\})\leq \frac{\bar{\mathbb{E}}[|X|^p]}{M^p}.
$$
\end{lemma}
\begin{theorem}\label{CTU}
Let $\{X^n\}_{n\in\mathbb{N}}\subset L^0(\Omega_T)$ be a sequence such that $X^n\uparrow X$, q.s., and there exists a $\mathbb{P}\in\mathcal{P}_G$, $E^\mathbb{P}[X^0]>-\infty$, then $\bar{\mathbb{E}}[X^n]\uparrow\bar{\mathbb{E}}[X]$.
\end{theorem}
\noindent Unlike the classical downwards monotone convergence theorem, the one in the $G$-framework only holds true for a sequence from a subset of $L^0(\Omega_T)$ (cf. Theorem 31 in Denis et al. \cite{DHP}).
\begin{theorem}\label{CTD}
Let $\{X^n\}_{n\in\mathbb{N}}\subset L_G^1(\Omega_T)$ be a sequence such that $X^n\downarrow X$, q.s., then $\bar{\mathbb{E}}[X^n]\downarrow\bar{\mathbb{E}}[X]$.
\end{theorem}
\noindent Moreover, by a classical argument, we have the following Fatou's lemma and the inequality of Jensen type in the $G$-framework. 
\begin{lemma}\label{le 2.2.12}
Assume that $\{X^n\}_{n\in\mathbb{N}}$ is a sequence in $L^0(\Omega_T)$ and for a $Y\in L^0(\Omega_T)$ that satisfies $\bar{\mathbb{E}}[|Y|]<+\infty$ and all $n\in\mathbb{N}$, $X^n \geq Y$, q.s., then
$$
\bar{\mathbb{E}}[\liminf_{n\rightarrow+\infty}X^n]\leq\liminf_{n\rightarrow+\infty}\bar{\mathbb{E}}[X^n].\\
$$
\end{lemma}
\noindent\textbf{Proof:} From (\ref{grep}) and by the classical Fatou-Lebesgue theorem,
we have for each $\mathbb{P}\in\mathcal{P}_G$, 
\begin{align*}
E^\mathbb{P}[\liminf_{n\rightarrow +\infty}X^{n}]\leq\liminf_{n\rightarrow +\infty}E^\mathbb{P}[X^{n}]
\leq\liminf_{n\rightarrow +\infty}\sup_{\mathbb{P}\in\mathcal{P}_G}E^\mathbb{P}[X^{n}]
=\liminf_{n\rightarrow +\infty}\bar{\mathbb{E}}[X^{n}].
\end{align*}
Taking the supremum of the left-hand side over all $\mathbb{P}\in\mathcal{P}_G$, we can easily obtain the desired result.\hfill{}$\square$
\begin{lemma}\label{le 2.2.2}Let $\rho : \mathbb{R} \rightarrow \mathbb{R}$ be an
increasing and concave function, then for each
$X\in L^0(\Omega_T),$ the following inequality holds:
$$
\bar{\mathbb{E}}[\rho(X)]\leq\rho(\bar{\mathbb{E}}[X]).$$
\end{lemma}
\noindent A representation theorem for $L^p_G(\Omega_T)$ can also be found in Denis et al. \cite{DHP}:
\begin{theorem}\label{lpg}
$$
L^p_G(\Omega_T)=\{X\in L^0(\Omega_T): X\ has\ a\ q.c.\ version,\ \lim_{N\rightarrow+\infty}\bar{\mathbb{E}}[|X|^p{\bf{1}}_{|X|>N}]=0\}.
$$
\end{theorem}
\noindent This dual definition of $L^p_G(\Omega_T)$ is more explicit to verify than its original definition given by the completion of $L^0_{ip}(\Omega_T)$.
\subsection{G-stochastic calculus}
In Peng \cite{P3}, generalized It\^o integrals with respect to $G$-Brownian motion and a generalized It\^o formula are established.
\begin{definition} \label{caocao}
A partition of $[0,T]$ is a finite
ordered subset $\pi^N_{[0, T]}=\{t_0, t_1,\ldots,t_N\}$ such that
$0=t_0<t_1<\ldots<t_N=T$. We set
$$\mu(\pi^N_{[0,T]}):=\max_{k=0, 1,\ldots,N-1}|t_{k+1}-t_k|.$$
For each $p\geq1$, we define
$$
M_G^{p,0}([0,T]):=\bigg\{ \eta_t=\sum_{k=0}^{N-1}\xi_{k}{\bf{1}}_{[t_{k},t_{k+1})}(t):
\xi_{k}\in L^0_{ip}(\Omega_{t_{k}})\bigg\},
$$
and we denote by $M_G^p([0,T])$ the completion of
$M_G^{p,0}([0,T])$ under the norm:
\begin{equation}\label{mnorm}\|\eta\|_{M_G^p([0,T])}:=\bigg(\frac{1}{T}\int_0^T\bar{\mathbb{E}}[|\eta_t|^p]dt\bigg)^\frac{1}{p}.
\end{equation}
\end{definition}
\begin{remark}\label{remconna}
By Definition \ref{caocao}, 
if $\eta$ is an element in $M^p_G([0,T])$, then there exists a sequence $\{\eta^n\}_{n\in\mathbb{N}}$ in $M_G^{p,0}([0,T])$, such that $
\lim\limits_{n\rightarrow+\infty}\int_0^T\bar{\mathbb{E}}[|\eta^n_t-\eta_t|^p]dt\rightarrow 0
$. It is readily observed that for $t\in[0,T]$, $\lambda$-a.e., 
$\bar{\mathbb{E}}[|\eta^n_t-\eta_t|^p]\rightarrow 0$ and thus, $\eta_t$ is an element in $L^p_G(\Omega_t)$, $\lambda$-a.e..
\end{remark}
\noindent Let $\mathbf{a}=(a_1,\ldots,a_d)^{\bf{Tr}}$ be a given vector in
$\mathbb{R}^d$ and
$B^\mathbf{a}=(\mathbf{a}, B),$ where
$(\mathbf{a},B)$ denotes the
scalar  product of $\mathbf{a}$ and  $B$.
\begin{definition}
For each $\eta\in M_G^{2,0}([0,T])$  with the form:
$$\eta_t=\sum_{k=0}^{N-1}\xi_{k}{\bf{1}}_{[t_{k},t_{k+1})}(t),$$ we
define
$$
\mathcal{I}_{[0, T]}(\eta)=\int_0^T\eta_tdB_t^\mathbf{a}
:=\sum_{k=0}^{N-1}\xi_k(B_{t_{k+1}}^\mathbf{a}-B_{t_k}^\mathbf{a}),
$$
and the mapping
can be continuously extended to $\mathcal{I}_{[0, T]}: M_G^2([0,T])\rightarrow L_G^2(\Omega_T)$.
Then, for each $\eta\in M_G^2([0,T])$, the stochastic integral is
defined by
$$\int_0^T\eta_tdB_t^\mathbf{a}:=\mathcal{I}_{[0, T]}(\eta).$$
\end{definition}
\noindent Let $\langle B^\mathbf{a}\rangle$ denote the
quadratic variation process of $B^\mathbf{a}$, which is formulated in $M^2_G([0,T])$ by
$$
\langle B^\mathbf{a}\rangle_t:=\lim_{\mu(\pi^N_{[0, T]})\rightarrow 0}\sum^{N-1}_{k=0}(B^\mathbf{a}_{t^N_{k+1}}-B^\mathbf{a}_{t^N_k})^2=(B^\mathbf{a}_t)^2-2\int^t_0B^\mathbf{a}_sdB^\mathbf{a}_s.
$$
We define
$$\sigma_{\mathbf{a}\mathbf{a}^{\bf{Tr}}}:=\sup_{\gamma\in\Gamma}tr(\gamma\gamma^{\bf{Tr}}\mathbf{a}\mathbf{a}^{\bf{Tr}}).$$
By Corollary 5.7 in Chapter III of Peng \cite{P3}, we have
$$\langle B\rangle_t\in t\Sigma:=\{t\times\gamma\gamma^{\bf{Tr}}:
\gamma\in\Gamma\},\ 0\leq t\leq T.$$ Therefore, for each $0\leq s\leq t\leq T,$
\begin{equation}\label{pquar}
\langle B^\mathbf{a}\rangle_t-\langle B^\mathbf{a}\rangle_s\leq \sigma_{\mathbf{a}\mathbf{a}^{\bf{Tr}}}(t-s).
\end{equation}
\begin{definition}
We define the mapping $\mathcal {Q}_{[0, T]}: M_G^{1,0}([0,T])\rightarrow L_G^1(\Omega_T)$ as follows:
$$\mathcal {Q}_{[0, T]}(\eta)=\int_0^T\eta_td\langle
B^\mathbf{a}\rangle_t:=\sum_{k=0}^{N-1}\xi_k(\langle
B^\mathbf{a}\rangle_{t_{k+1}}-\langle
B^\mathbf{a}\rangle_{t_{k}}),$$ and we extend it to $\mathcal {Q}_{[0, T]}: M_G^{1}([0,T])\rightarrow L_G^1(\Omega_T)$. This extended mapping 
defines 
$\int_0^T\eta_sd\langle B^\mathbf{a}\rangle_s$ for each $\eta\in M_G^1([0,T]).$
\end{definition}
\noindent  For two given vectors $\mathbf{a}$, $\bar{\mathbf{a}}\in\mathbb{R}^d$, the mutual variation process of $B^\mathbf{a}$ and
$B^{\bar{\mathbf{a}}}$ is defined by $$\langle
B^\mathbf{a},B^{\bar{\mathbf{a}}} \rangle_t:=\frac{1}{4}(\langle
B^{\mathbf{a}+\bar{\mathbf{a}}}\rangle_t-\langle
B^{\mathbf{a}-\bar{\mathbf{a}}}\rangle_t).$$
Then, for each $\eta\in M^1_G([0, T])$,
$$
\int_0^T\eta_td\langle
B^\mathbf{a},B^{\bar{\mathbf{a}}}\rangle_t:=\frac{1}{4}\bigg(\int_0^T\eta_td\langle
B^{\mathbf{a}+\bar{\mathbf{a}}}\rangle_t-\int_0^T\eta_td\langle
B^{\mathbf{a}-\bar{\mathbf{a}}}\rangle_t\bigg).
$$
\noindent In view of the dual formulation of $G$-expectation (\ref{grep}) and the property of the quadratic variation process $\langle B\rangle$ in the $G$-framework (\ref{pquar}), the following BDG type inequalities are obvious. (cf. Theorem 2.1 and 2.2 in Gao \cite{G}) 
\begin{lemma}\label{th 2.1.2}
Let $p\geq1$, $\mathbf{a},\bar{\mathbf{a}}\in\mathbb{R}^d$, $\eta\in M_G^p([0,T])$ and $0\leq s\leq t\leq T$. Then, 
\begin{align*}
\bar{\mathbb{E}}\bigg[\sup_{s\leq u\leq
t}\bigg|\int^u_s\eta_rd\langle B^{\mathbf{a}}, B^{\bar{\mathbf{a}}}\rangle_r\bigg|^p\bigg]\leq
\bigg(\frac{\sigma_{(\mathbf{a}+\bar{\mathbf{a}})(\mathbf{a}+\bar{\mathbf{a}})^{\bf{Tr}}}+\sigma_{(\mathbf{a}-\bar{\mathbf{a}})(\mathbf{a}-\bar{\mathbf{a}})^{\bf{Tr}}}}{4}\bigg)^p
(t-s)^{p-1}\int_s^t\bar{\mathbb{E}}[|\eta_u|^p]du.
\end{align*}
\end{lemma}
\begin{lemma}\label{th 2.1.1}Let $p\geq2$, $\mathbf{a}\in\mathbb{R}^d$, $\eta\in M_G^p([0,T])$ and $0\leq s\leq t\leq T$. Then,
$$
\bar{\mathbb{E}}\bigg[\sup_{s\leq u\leq
t}\bigg|\int^u_s\eta_rdB^{\mathbf{a}}_r\bigg|^p\bigg]\leq
C_p\sigma^{p/2}_{\mathbf{a}\mathbf{a}^{\bf{Tr}}}|t-s|^{\frac{p}{2}-1}\bigg(\int_s^t\bar{\mathbb{E}}[|\eta_u|^p]du\bigg),
$$
where $C_p>0$ is a constant independent of $\mathbf{a}, \
\eta$ and $\Gamma.$
\end{lemma}
\noindent At the end of this subsection, we introduce the following $G$-It\^o formula that can be found as Proposition 6.3 in Chapter III of Peng \cite{P3}. For each $0\leq s\leq t\leq T$, consider an $n$-dimensional $G$-It\^o process:
$$
X^{\nu}_t=X^{\nu}_s+\int^t_s b^{\nu}_udu+\sum^d_{i, j=1}\int^t_s h^{\nu ij}_u d\langle B^i,B^j\rangle_u+\sum^d_{j=1}\int^t_s g^{\nu j}_u dB^j_u,\ \nu= 1,\ldots,n.
$$
\begin{lemma}
Let $\Phi\in\mathcal{C}^2(\mathbb{R}^n)$ be a real function with bounded derivatives such that $\{\partial^2_{x^\mu x^\nu}\Phi\}^n_{\mu, \nu=1}$ are uniformly Lipschitz. Let $b^{\nu}$, $h^{\nu ij}$ and $g^{\nu j}\in M^2_G([0,T])$, $\nu= 1,\ldots,n$, $i, j= 1,\ldots,d,$ be bounded processes. Then, we have
\begin{align}\label{itofu}
\Phi(X_t)-\Phi(X_s)&=\int^t_s\partial_{x^\nu}\Phi(X_u)b^{\nu}_{u}du
+\int^t_s\partial_{x^\nu}\Phi(X_u)h^{\nu ij}_{u}d\langle B^i, B^j\rangle_u\\
&+\int^t_s\partial_{x^\nu}\Phi(X_u)g^{\nu j}_{u}dB^j_u+\frac{1}{2}\int^t_s\partial^2_{x^\mu x^\nu}\Phi(X_u)g^{\mu i}_u g^{\nu j}_ud\langle B^i, B^j\rangle_u,\notag
\end{align}
in which the equality holds in the sense of $L^2_G(\Omega_t)$.
\end{lemma}
\begin{remark}
In (\ref{itofu}), we adopt the Einstein convention, i.e., the repeated indices $\nu$, $\mu$, $i$ and $j$ imply the summation. \end{remark}
\section{Solvability of GSDEs with integral-Lipschitz coefficients}
\noindent In this section, we give our main result of this paper, that is, the
existence and uniqueness theorems for GSDEs with integral-Lipschitz coefficients. From now on, $C$ denotes a positive constant whose value may vary from line to line.
\subsection{Formulation to GSDEs and assumptions} We rewrite (\ref{eq 2}) into the following form: 
\begin{align}\label{eq 12}
X(t)=x+\int^t_0 b(s,X(s))ds +\sum^d_{i,j=1}\int^t_0 h_{ij}(s&,X(s))d\langle B^i, B^j\rangle_s\\
&+ \sum^d_{j=1}\int^t_0 g _j(s,X(s)) dB^j_s,\ 0\leq t\leq T,\notag
\end{align}
where the initial value $x\in \mathbb {R}^n$ is a given vector, $B^i$ is the $i$th component of a $d$-dimensional $G$-Brownian motion, $\langle B^i, B^j\rangle$ is the mutual variation process of $B^i$ and $B^j$ and $b$, $h_{ij}$, $g_j$ are given functions that satisfy for each $x\in \mathbb{R}^n$, $b(\cdot, x)$, $h_{ij}(\cdot, x)$, $g_j(\cdot, x)\in M^2_G([0,T];\mathbb{R}^n)$, $i$, $j=1, \ldots, d$.\\[6pt]
We now state our main assumptions on the coefficients of GSDE (\ref{eq 12}), which will be our main interest in the sequel:
\begin{assumption}\label{ass1}
 For each $t\in[0,T]$ and $x$, $x_1$, $x_2\in \mathbb{R}^n$,
\begin{mnlist}
\item\label{H1}$|b(t,x_1)-b(t,x_2)|^2+|h(t,x_1)-h(t,x_2)|^2+|g(t,x_1)-g(t,x_2)|^2\leq|\beta(t)|^2\rho(|x_1-x_2|^2)$;
\item\label{H2}$|b(t,x)|^2+|h(t,x)|^2+|g(t,x)|^2\leq|\beta_1(t)|^2+\beta^2_2|x|^2$,
\end{mnlist}
where $\beta:[0,T]\rightarrow\mathbb{R}^+$ is square integrable, $\beta_1\in M_G^2([0,T])$, $\beta_2\in\mathbb{R}^+$ and  $\rho:(0,+\infty)\rightarrow(0,+\infty)$ is a continuous
increasing and concave function that vanishes at $0+$ and satisfies 
\begin{equation}\label{eq 11}
\int_0^1\frac{dr}{\rho(r)}=+\infty.
\end{equation}
\end{assumption}
\begin{assumption}\label{ass2}
For each $t\in[0,T]$ and $x$, $x_1$, $x_2\in \mathbb{R}^n$,
\begin{mnnlist}
\item \label{H1'}$\left\{\begin{array}{l}
|b(t,x_1)-b(t,x_2)|+|h(t,x_1)-h(t,x_2)|\leq\beta(t)\rho_1(|x_1-x_2|);\\[3pt]
|g(t,x_1)-g(t,x_2)|^2\leq|\beta(t)|^2\rho_2(|x_1-x_2|^2);
\end{array}\right.$
\item \label{H2'}$|b(t,x)|^p+|h(t,x)|^p+|g(t,x)|^p\leq
|\beta_1(t)|^p+\beta^p_2|x|^p$,
\end{mnnlist}
where $\beta:[0,T]\rightarrow\mathbb{R}^+$ is square integrable, for some $p>2$, $\beta_1\in M^p_G([0,T])$, $\beta_2\in \mathbb{R}^+$ and both $\rho_1$, $\rho_2:(0,+\infty)\rightarrow(0,+\infty)$ are continuous increasing and concave functions that vanish at $0+$ and satisfy (\ref{eq 11}).
We assume moreover that $$\rho_3(r):=\frac{\rho_2(r^2)}{r},\ r\in(0,+\infty),$$ is also a continuous increasing and concave function that vanishes at $0+$ and satisfies
$$\int_0^1\frac{dr}{\rho_1(r)+\rho_3(r)}=+\infty.$$
\end{assumption}
\begin{remark}We give an example to show that (H\ref{H1'}') is ``weaker'' than (H\ref{H1}). If we set
\begin{align*}
\left\{\begin{array}{l}
\rho_1(r)=r\ln\frac{1}{r};\\[3pt]
\rho_2(r)=r\ln\frac{1}{r},
\end{array}\right.\ r\in(0,+\infty),
\end{align*}
then (H\ref{H1}') is satisfied but (H\ref{H1'}) is not.
\end{remark}
\noindent To ensure that (\ref{eq 12}) is well defined, all the integrands in (\ref{eq 12}) should be in $M^2_G([0, T]; \mathbb{R}^n)$. Thus, we need the following lemma:
\begin{lemma}\label{complex}
For some $q\geq 1$, $\zeta$ is a function that satisfies  for each $x\in \mathbb{R}^n$, $\zeta(\cdot, x)\in M^q_G([0,T];\mathbb{R}^n)$. We assume moreover that, for each $x$, $x_1$, $x_2\in \mathbb{R}^n$:
\begin{nmnnlist}
\item\label{A1} $|\zeta(t,x_1)-\zeta(t,x_2)|\leq \beta(t)\gamma(|x_1-x_2|)$;
\item\label{A2} $|\zeta(t,x)|\leq |\beta_1(t)|+\beta_2|x|$,
\end{nmnnlist}
where $\beta: [0,T]\rightarrow\mathbb{R}^+$ is $q$-integrable, $\beta_1\in M^q_G([0,T])$, $\beta_2\in\mathbb{R}^+$ and $\gamma:(0, +\infty)\rightarrow (0,+\infty)$ is an increasing function vanishes at $0$. Then, for each $X\in M^q_G([0,T];\mathbb{R}^n)$, $\zeta(\cdot, X_\cdot)$ is an element in $M^q_G([0,T];\mathbb{R}^n)$.
\end{lemma}
\begin{remark}When $q=2$,
all the coefficients in GSDE (\ref{eq 12}) satisfy both (A\ref{A1}) and (A\ref{A2}) under either Assumption \ref{ass1} or \ref{ass2}.
Therefore, 
the $G$-stochastic integrals in GSDE (\ref{eq 12}) are well defined for any solution $X\in M^2_G([0,T];\mathbb{R}^n)$. We postpone the proof of this lemma to the appendix. 
\end{remark}
\subsection{Main result}
As a starting point, we first refer to an inequality in Bihari \cite{B} (Bihari's inequality). Then, we prove the existence and uniqueness theorem for the GSDE (\ref{eq 12}) under Assumption \ref{ass1}.
\begin{lemma} \label{le 2.2.1}
Let $\rho:(0,+\infty)\rightarrow(0,+\infty)$ be a continuous and
increasing function that vanishes at $0+$ and satisfies (\ref{eq 11}).
Let $u$ be a measurable and non-negative function defined on $(0,+\infty)$ that satisfies 
$$
u(t)\leq a+\int^t_0 \kappa(s)\rho(u(s))ds,\ t\in(0,+\infty),
$$
where $a\in \mathbb{R}^+$ and $\kappa: [0,T]\rightarrow\mathbb{R}^+$ is Lebesgue
integrable. We have
\begin{mmmylist}
\item If $a=0$, then $u(t)=0$, $t\in(0, +\infty)$, $\lambda$-a.e.;
\item If $a>0$, we define
$$v(t):=\int^t_{t_0} \frac{1}{\rho(s)}ds,\ t\in\mathbb{R}^+,$$
where $t_0\in (0,+\infty)$, then
\begin{equation*}
u(t)\leq v^{-1} \bigg(v(a)+\int^t_0 \kappa(s) ds\bigg).
\end{equation*}
\end{mmmylist}
\end{lemma}
\begin{theorem}\label{th 3.1}
Under Assumption \ref{ass1} there exists a unique process $X\in M^2_G([0,T];\mathbb{R}^n)$ that satisfies the GSDE (\ref{eq 12}).
\end{theorem}
\noindent\textbf{Proof:} We begin with the proof of the uniqueness. Suppose $X(\cdot; x_i)\in M^2_G([0,T];\mathbb{R}^n)$ is a solution to the GSDE (\ref{eq 12}) with initial value $x_i$, $i=1, 2$, then we calculate
\begin{align*}
|X(t;x_1)-X(t;x_2)|^2&\leq C\bigg(|x_1-x_2|^2+\bigg|\int_0^t(b(s,X(s;x_1))-b(s,X(s;x_2)))ds\bigg|^2\\
&+\bigg|\sum^d_{i,j=1}\int_0^t(h_{ij}(s,X(s;x_1))-h_{ij}(s,X(s;x_2)))d\langle B^i,B^j\rangle_s\bigg|^2\\
&+\bigg|\sum^d_{j=1}\int_0^t(g_j(s,X(s;x_1))-g_j(s,X(s;x_2)))dB^j_s\bigg|^2\bigg).
\end{align*}
By the BDG type inequalities and (H\ref{H1}), we deduce
\begin{align*}
\bar{\mathbb{E}}\bigg[\sup_{0\leq s\leq t}\bigg|\int_0^s
(b(r,X(r;x_1))&-b(r,X(r;x_2)))dr\bigg|^2\bigg]\\
&\leq C\int_0^t|\beta(s)|^2\bar{\mathbb{E}}[\rho(|X(s;x_1)-X(s;x_2)|^2)]ds;
\end{align*}
\begin{align*}
\bar{\mathbb{E}}\bigg[\sup_{0\leq s\leq t}\bigg|\int_0^s
(h_{ij}(r,X(r;x_1))&-h_{ij}(r,X(r;x_2)))d\langle B^i,B^j\rangle_r\bigg|^2\bigg]\\&\leq C\int_0^t|\beta(s)|^2\bar{\mathbb{E}}[\rho(|X(s;x_1)-X(s;x_2)|^2)]ds;
\end{align*}
and
\begin{align*}
\bar{\mathbb{E}}\bigg[\sup_{0\leq s\leq t}\bigg|\int_0^s
(g_j(r,X(r;x_1))&-g_j(r,X(r;x_2)))dB^j_r\bigg|^2\bigg]\\
&\leq C\int_0^t|\beta(s)|^2\bar{\mathbb{E}}[\rho(|X(s;x_1)-X(s;x_2)|^2)]ds.
\end{align*}
Set
$$
u(t):=\sup_{0\leq s\leq t}\bar{\mathbb{E}}[|X(s;x_1)-X(s;x_2)|^2],
$$
then
$$u(t)\leq C\bigg(|x_1-x_2|^2+\int_0^t|\beta(s)|^2\bar{\mathbb{E}}[\rho(|X(s;x_1)-X(s;x_2)|^2)]ds\bigg).$$
As $\rho$ is an increasing and concave function, by Lemma \ref{le
2.2.2}, we have
\begin{align*}
u(t)&\leq C\bigg(|x_1-x_2|^2+\int_0^t|\beta(s)|^2\rho(\bar{\mathbb{E}}[|X(s;x_1)-X(s;x_2)|^2])ds\bigg)\\
&\leq C\bigg(|x_1-x_2|^2+\int_0^t|\beta(s)|^2\rho(\sup_{0\leq r\leq s}\bar{\mathbb{E}}[|X(r;x_1)-X(r;x_2)|^2])ds\bigg)\\
&\leq C\bigg(|x_1-x_2|^2+\int_0^t|\beta(s)|^2\rho(u(s))ds\bigg).
\end{align*}
By Lemma \ref{le 2.2.1}, we obtain
$$
u(t)\leq v^{-1} \bigg(v(C|x_1-x_2|^2)+C\int^t_0|\beta(s)|^2ds\bigg).
$$
In particular, if $x_1=x_2,$ then $u(t)=0$, $0\leq t\leq T$, which implies the pathwise uniqueness.\\[6pt]
Now we start to prove the existence. We
define a Picard sequence $\{X^m(\cdot)\}_{m\in\mathbb{N}}$ by the following procedure:
$$X^0(t)=x,\ 0\leq t\leq T;$$
and
\begin{align}\label{eq m+1}
X^{m+1}(t)=x&+\int_0^tb(s,X^m(s))ds\\
&+\sum^d_{i,j=1}\int_0^th_{ij}(s,X^m(s))d\langle
B^i,B^j\rangle_s
+\sum^d_{j=1}\int_0^tg_j(s,X^m(s))dB^j_s,\ 0\leq t\leq T.\notag
\end{align}
By Lemma \ref{complex}, the sequence
$\{X^m(\cdot)\}_{m\in\mathbb{N}}$ is well
defined in $M_G^2([0,T];\mathbb{R}^n)$.\\[6pt]
First, we establish an a priori estimate for
$\{\bar{\mathbb{E}}[|X^m(t)|^2]\}_{m\in\mathbb{N}}$.
From (\ref{eq m+1}), by the BDG type inequalities, we deduce
\begin{align*}
\bar{\mathbb{E}}[|X^{m+1}(t)|^2]&\leq
C\bigg(|x|^2+\int_0^t\bar{\mathbb{E}}[|\beta_1(s)|^2+\beta_2^2|X^m(s)|^2]ds\bigg)\\
&\leq C\bigg(|x|^2+\int_0^t\bar{\mathbb{E}}[|\beta_1(s)|^2]ds+\beta_2^2\int_0^t\bar{\mathbb{E}}[|X^m(s)|^2]ds\bigg).
\end{align*}
Set
\begin{equation*}\label{eq 13}
p(t) :=Ce^{C\beta^2_2t}\bigg(|x|^2+\int_0^t\bar{\mathbb{E}}[|\beta_1(s)|^2]ds\bigg),
\end{equation*}
then $p(\cdot)$ is the solution to the following ordinary differential equation:
\begin{equation*}\label{eq 14}
p(t)=C\bigg(|x|^2+\int_0^t\bar{\mathbb{E}}[|\beta_1(s)|^2]ds+\beta_2^2\int_0^tp(s)ds\bigg).
\end{equation*}
By recurrence, it is easy to verify that for each $m\in\mathbb{N}$,
$$\bar{\mathbb{E}}[|X^m(t)|^2]\leq p(t),$$
the right-hand side of which is continuous and thus, bounded on $[0, T]$.\\[6pt]
Secondly, for each $k$, $m\in\mathbb{N}$, we define
$$u_{k+1,m}(t):=\sup_{0\leq s\leq t}\bar{\mathbb{E}}[|X^{k+1+m}(s)-X^{k+1}(s)|^2].$$
By the definition of the sequence $\{X^m(\cdot)\}_{m\in\mathbb{N}}$, we have
\begin{align*}
X^{k+1+m}(t)-X^{k+1}(t)&=\int_0^t(b(s,X^{k+m}(s))-b(s,X^k(s)))ds\\
&+\sum^d_{i,j=1}\int_0^t(h_{ij}(s,X^{k+m}(s))-h_{ij}(s,X^k(s)))d\langle B^i,B^j\rangle_s\\
&+\sum^d_{j=1}\int_0^t(g_j(s,X^{k+m}(s))-g_j(s,X^k(s)))dB^j_s.
\end{align*}
By an argument similar to the one in the proof of the uniqueness, we obtain
$$u_{k+1,m}(t)\leq C\int_0^t|\beta(s)|^2\rho(u_{k,m}(s))ds.$$
Set
\begin{equation*}
v_k(t):=\sup_{m\in\mathbb{N}}u_{k,m}(t),\ 0\leq t\leq T,
\end{equation*}
then
\begin{equation}\label{new1}0\leq v_{k+1}(t)\leq C\int_0^t|\beta(s)|^2\rho(v_{k}(s))ds.\end{equation}
Finally, we define
$$\alpha(t):=\limsup_{k\rightarrow+\infty}v_k(t),\ 0\leq
t\leq T,$$
which is uniformly bounded by $4p(t)$. Applying 
the Fatou-Lebesgue theorem to (\ref{new1}), we have
$$0\leq\alpha(t)\leq C\int_0^t\beta^2(s)\rho(\alpha(s))ds.
$$
By Lemma \ref{le 2.2.1}, we deduce
$$\alpha(t)=0,\ 0\leq t\leq T,$$ which implies that 
$\{X^m(\cdot)\}_{m\in\mathbb{N}}$ is a Cauchy sequence under the norm  
$\sup\limits_{0\leq t\leq T}(\bar{\mathbb{E}}[|\cdot|^2])^{\frac{1}{2}}$, which is stronger than the $M^2_G([0, T]; \mathbb{R}^n)$ norm (\ref{mnorm}). Therefore, one can find a process $X\in M_G^2([0,T];\mathbb{R}^n)$ that satisfies 
$$
\sup\limits_{0\leq t\leq T}\bar{\mathbb{E}}[|X^{m}(t)-X(t)|^2]\rightarrow 0, \ {\rm{as}}\ m\rightarrow +\infty.
$$
\noindent Moreover, it is readily observed that
\begin{align}\label{new 34}
\notag&\ \bar{\mathbb{E}}\bigg[\sup_{0\leq t\leq
T}\bigg|\int_0^t(b(s,X^{m}(s))-b(s,X(s)))ds\bigg|^2\bigg]\\\notag+& 
\sum_{i, j=1}^d\bar{\mathbb{E}}\bigg[\sup_{0\leq t\leq
T}\bigg|\int_0^t(h_{ij}(s,X^{m}(s))-h_{ij}(s,X(s)))d\langle B^i,B^j\rangle_s\bigg|^2\bigg]\\
+&\ \sum_{i=1}^d\bar{\mathbb{E}}\bigg[\sup_{0\leq t
\leq T}\bigg|\int_0^t(g_j(s,X^{m}(s))-g_j(s,X(s)))dB^j_s\bigg|^2\bigg]\\\notag
\leq &\ C\int_0^T|\beta(t)|^2\rho(\bar{\mathbb{E}}[|X^{m}(t)-X(t)|^2])dt\\[6pt]\notag
\leq &\ C\rho(\sup_{0\leq t\leq T}\bar{\mathbb{E}}[|X^{m}(t)-X(t)|^2]).
\end{align}
By the continuity of $\rho$ and $\rho(0+)=0$, we know that $\rho(\sup\limits_{0\leq t\leq T}\bar{\mathbb{E}}[|X^{m}(t)-X(t)|^2])\rightarrow 0$ and the left-hand side of (\ref{new 34}) converges to $0$. 
Thus, $\{X^m(\cdot)\}_{m\in\mathbb{N}}$ is a successive approximation to $X$, which is a solution to the GSDE (\ref{eq 12}) in $M_G^2([0,T];\mathbb{R}^n)$.\hfill{} $\square$\\[6pt]
In what follows, we give the existence and uniqueness theorem to the GSDE (\ref{eq 12}) under Assumption \ref{ass2} instead of Assumption \ref{ass1}. 
\begin{theorem}\label{th 3.5}
Under Assumption \ref{ass2} there exists a unique process $X\in M^2_G([0,T];\mathbb{R}^n)$ that satisfies GSDE (\ref{eq 12}).
\end{theorem}
\noindent\textbf{Proof:} We start with the proof of existence. Similar to (\ref{eq m+1}), we define a sequence of processes $\{X^m\}_{m\in\mathbb{N}}$ as
follows:
$$X^0(t)=x,\ 0\leq t\leq T;$$
and
\begin{align}\label{Picard2}
X^{m+1}(t)=x+\int_0^tb(s,X^m(s))ds
&+\sum^d_{i,j=1}\int_0^th_{ij}(s,X^m(s))d\langle
B^i, B^j\rangle_s\\
&+\sum^d_{j=1}\int_0^tg_j(s,X^{m+1}(s))d B^j_s,\ 0\leq t\leq T.\notag
\end{align}
Thanks to Theorem \ref{th 3.1}, the sequence $\{X^m\}_{m\in\mathbb{N}}$ is well defined in $M^2_G([0,T]; \mathbb{R}^n)$. \\[6pt]
We notice that the coefficients in (\ref{Picard2}) could not be bounded. In order to apply the $G$-It\^{o} formula, we shall firstly construct, for each $m\in\mathbb{N}$, a sequence of $G$-It\^o processes that approximates $X^m$, and whose coefficients are all truncated. These sequences are given by the following steps:\\[6pt]
\noindent Step 1: 
For each $N\in\mathbb{N}$, we set
\begin{align}\label{trunccoe}
\zeta^N(t, x)=
\left\{\begin{array}{r@{,}l}
\zeta(t, x)&\ {\rm{if}}\ |\zeta(t, x)|\leq N;\\[3pt]
\frac{N\zeta(t, x)}{|\zeta(t, x)|}&\ {\rm{if}}\ |\zeta(t, x)|> N,
\end{array}\right.
\end{align}
where $\zeta=b$, $h_{ij}$ or $g_j$, $i$, $j=1, \ldots, d$, respectively.
It is easy to verify that $b^N$, $h_{ij}^N$ and $g^N_j$ still satisfy (H\ref{H1'}') and (H\ref{H2'}'). \\[6pt]
Step 2: For each $m\in\mathbb{N}$, we define  
\begin{align*}
X^{m+1, N}(t)=x&+\int_0^tb^N(s,X^m(s))ds\\
&+\sum^d_{i,j=1}\int_0^th^N_{ij}(s,X^m(s))d\langle
B^i, B^j\rangle_s
+\sum^d_{j=1}\int_0^tg^N_j(s,X^{m+1}(s))d B^j_s,\ 0\leq t\leq T.
\end{align*}
By Lemma \ref{complex}, the sequence $\{X^{m, N}(\cdot)\}_{N\in\mathbb{N}}$ is also well defined in $M^2_G([0,T]; \mathbb{R}^n)$.\\[6pt]
Let us now establish an a priori estimate for $\{\bar{\mathbb{E}}[|X^{m}(t)|^p]\}_{m\in\mathbb{N}}$.
By (H{\ref{H2'}'}) and the BDG type inequalities,
\begin{align}\label{est2}
\bar{\mathbb{E}}[|X^{m+1}(t)|^p]\leq C\bigg(|x|^p&+\int_0^t\bar{\mathbb{E}}[|\beta_1(s)|^p]ds\\
&+\beta^p_2\int_0^t\bar{\mathbb{E}}[|X^{m}(s)|^p]ds+\beta^p_2\int_0^t\bar{\mathbb{E}}[|X^{m+1}(s)|^p]ds\bigg).\notag
\end{align}
By induction, we obtain that $\bar{\mathbb{E}}[|X^m(t)|^p]\leq p'(t)$, where $p'(\cdot)$ is the solution to the following ordinary differential equation:
\begin{equation*}
p'(t)=C\bigg(|x|^p+\int_0^t\bar{\mathbb{E}}[|\beta_1(s)|^p]ds+\beta_2^p\int_0^t p'(s)ds\bigg).
\end{equation*}
Since $p'(\cdot)$ is continuous and bounded on $[0, T]$, we have
\begin{align}\label{est1}
\sup_{m\in\mathbb{N}}\sup_{0\leq t\leq T}\bar{\mathbb{E}}[|X^m(t)|^p]\leq M <+\infty.
\end{align}
Fixing an $m>0$, we calculate
\begin{align*}
\sup_{0\leq t\leq T}\bar{\mathbb{E}}[|X^{m, N}(t)-X^m(t)|]
&\leq \bar{\mathbb{E}}\bigg[\int_0^T|b^N(t,X^m(t))-b(t,X^m(t))|dt\bigg]\\
&+\sum^d_{i,j=1}\bar{\mathbb{E}}\bigg[\int_0^T|h^N_{ij}(t,X^m(t))-h_{ij}(t,X^m(t))|d\langle
B^i, B^j\rangle_t\bigg]\\
&+\sup_{0\leq t\leq T}\sum^d_{j=1}\bar{\mathbb{E}}\bigg[\bigg|\int_0^t(g^N_j(s,X^{m+1}(s))-g_j(s,X^{m+1}(s)))d B^j_s\bigg|\bigg].
\end{align*}
By the definition of the truncated coefficients and the BDG type inequalities, we deduce
\begin{align}\label{est10}
\sup_{0\leq t\leq T}\bar{\mathbb{E}}[|X^{m, N}(t)-X^m(t)|]\notag
&\leq\int_0^T\bar{\mathbb{E}}[|b(t,X^m(t))|{\bf{1}}_{\{|b(t, X^m(t))|> N\}}]dt\\
&+C\bigg(\sum^d_{i, j=1}\int_0^T\bar{\mathbb{E}}[|h_{ij}(t,X^m(t))|{\bf{1}}_{\{|h_{ij}(t, X^m(t))|> N\}}]dt\\
&+\sum^d_{j=1}\bigg(\int_0^T\bar{\mathbb{E}}[|g_{j}(t,X^{m+1}(t))|^2{\bf{1}}_{\{|g_{j}(t, X^{m+1}(t))|> N\}}]dt\bigg)^{\frac{1}{2}}\bigg).\notag
\end{align}
By Lemma \ref{complex}, for each $m\in\mathbb{N}$, $b(\cdot, X^m_\cdot)$, $h_{ij}(\cdot, X^m_\cdot)$, $g_j(\cdot, X^m_\cdot)\in M^2_G([0, T]; \mathbb{R}^n)$, $i$, $j=1, \ldots, d$. Then, by Remark \ref{remconna} and Theorem \ref{lpg} along with Lebesgue's dominated convergence theorem, the right-hand side of (\ref{est10}) converges to 0. Therefore, 
\begin{equation}\label{con1}\sup_{0\leq t\leq T}\bar{\mathbb{E}}[|X^{m, N}(t)-X^m(t)|]\rightarrow 0,\ {\rm{as}}\ N\rightarrow +\infty.
\end{equation}
Since $|x|$ is not a $\mathcal{C}^2(\mathbb{R}^n)$ function, we have to approximate $|x|$ by a sequence of $\mathcal{C}^2(\mathbb{R}^n)$ functions, i.e., 
$\{F_\varepsilon(x)\}_{\varepsilon>0}$, where
$$F_\varepsilon(x):=(|x|^2+\varepsilon)^\frac{1}{2},\ x\in\mathbb{R}^n.$$
We notice that
\begin{equation}\label{prop1}
F_\varepsilon(x)\geq \varepsilon^{\frac{1}{2}};\ \bigg|\frac{\partial F_\varepsilon(x)}{\partial x_i}\bigg|\leq1;\
\bigg|\frac{\partial^2F_\varepsilon(x)}{\partial x_i \partial x_j}\bigg|\leq\frac{2}{F_\varepsilon(x)};
\end{equation}
and thus, $\frac{\partial F_\varepsilon(x)}{\partial x_i}$, $\frac{\partial^2F_\varepsilon(x)}{\partial x_i \partial x_j}$, $i, j=1, \ldots, n$, are uniformly Lipschitz.\\[6pt]
Fixing an $\varepsilon\in(0, +\infty)$, we define
$$\Delta F_\varepsilon^{k, m, N}(t):= F_\varepsilon(\Delta X^{k, m, N}(t))-F_\varepsilon(\Delta X^{k, m}(t)), $$
where
$$\Delta X^{k, m, N}(t)= X^{k+m, N}(t)-X^{k, N}(t)$$
and
$$\Delta X^{k, m}(t)=X^{k+m}(t)-X^{k}(t).$$
We apply
the $G$-It\^o formula to $F_\varepsilon(\Delta X^{k+1, m, N}(t))$ and take $G$-expectation on both sides. Then, from (\ref{prop1}) and by the BDG type inequalities, it is easy to show that
\begin{align}\label{eq key}
\bar{\mathbb{E}}[F_\varepsilon(\Delta X^{k+1, m, N}(t))]
&\leq\int_0^t\bar{\mathbb{E}}[|b^N(s,X^{k+m}(s))-b^N(s,X^{k}(s))|]ds\notag\\
&+C\bigg(\sum^d_{i,j=1}\int_0^t\bar{\mathbb{E}}[|h^N_{ij}(s,X^{k+m}(s))-h^N_{ij}(s,X^{k}(s))|]ds\\
&+\sum^d_{j=1}\int_0^t\bar{\mathbb{E}}\bigg[\frac{|g^N_j(s,X^{k+m+1}(s))-g^N_j(s,X^{k+1}(s))|^2}{F_\varepsilon(\Delta X^{k+1, m, N}(s))}\bigg]ds\bigg),\ 0\leq t\leq T.\notag
\end{align}
By (H\ref{H1'}') and Lemma \ref{le 2.2.2}, we have 
\begin{align}\label{new89}
\bar{\mathbb{E}}[F_\varepsilon(\Delta X^{k+1, m, N}(t))]\leq C\int_0^t\beta(s)\bigg(\rho_1(\bar{\mathbb{E}}[|\Delta X^{k, m}(s)|])+\bar{\mathbb{E}}\bigg[\frac{\rho_2(|\Delta X^{k+1, m}(s)|^2)}{F_\varepsilon(\Delta X^{k+1, m, N}(s))}\bigg]\bigg)ds.
\end{align}
From (\ref{prop1}), we know that $F_\varepsilon(x)$ is uniformly Lipschitz. Based on this fact and (\ref{con1}), we obtain
\begin{align}\label{est7}
\sup_{0 \leq t\leq T}\bar{\mathbb{E}}[|\Delta F_\varepsilon^{k+1, m, N}(t)|]&\leq\sup_{0 \leq t\leq T}\bar{\mathbb{E}}
[|\Delta X^{k+1, m, N}(t)-\Delta X^{k+1, m}(t)|]\notag\\
&\leq\sup_{0 \leq t\leq T}\bar{\mathbb{E}}
[|X^{k+1+m, N}(t)-X^{k+1+m}(t)|]\\&+\sup_{0 \leq t\leq T}\bar{\mathbb{E}}
[|X^{k+1, N}(t)-X^{k+1}(t)|]\rightarrow 0,\ {\rm{as}}\ N\rightarrow +\infty.\notag
\end{align}
Since $\rho_1$, $\rho_2: (0, +\infty) \rightarrow (0, +\infty)$ are concave and vanish at 0+, for each $\delta\in(0, +\infty)$, we can find a positive constant $K_\delta$ such that for each $x\in[\delta, +\infty)$, $\rho_1(x)$, $\rho_2(x)\leq K_\delta x$. Fixing a $\delta>0$ and $M\in(\delta, +\infty)$, we calculate
\begin{align}\label{est18}
\sup_{0\leq t\leq T}\bar{\mathbb{E}}\bigg[
\bigg|\frac{\rho_2(|\Delta X^{k+1,m}(t)|^2)}{F_\varepsilon(\Delta X^{k+1, m, N}(t))}&-\frac{\rho_2(|\Delta X^{k+1,m}(t)|^2)}{F_\varepsilon(\Delta X^{k+1, m}(t))}\bigg|\bigg]\notag\\
&\leq 2\varepsilon^{-\frac{1}{2}} K_\delta\sup_{0\leq t\leq T}\bar{\mathbb{E}}[|\Delta X^{k+1,m}(t)|^2{\bf{1}}_{\{|\Delta X^{k+1,m}(t)|^2> M\}}]\\
&+\varepsilon^{-1}\rho_2(M)\sup_{0 \leq t\leq T}\bar{\mathbb{E}}[|\Delta F_\varepsilon^{k+1, m, N}(t)|].\notag
\end{align}
On account of (\ref{est7}) and by H\"older's inequality and  Lemma \ref{markov}, 
\begin{align*}
\limsup_{N\rightarrow+\infty}\sup_{0\leq t\leq T}\bar{\mathbb{E}}\bigg[\bigg|\frac{\rho_2(|\Delta X^{k+1, m}(t)|^2)}{F_\varepsilon(\Delta X^{k+1, m, N}(t))}&-\frac{\rho_2(|\Delta X^{k+1, m}(t)|^2)}{F_\varepsilon(\Delta X^{k+1, m}(t))}\bigg|\bigg]\\
&\leq\frac{2K_\delta}{\varepsilon^{\frac{1}{2}}{M^{p-2}}}\sup_{0\leq t\leq T}\bar{\mathbb{E}}[
|\Delta X^{k+1, m}(t)|^p].
\end{align*}
As $M$ can be arbitrary large and $\bar{\mathbb{E}}[
|\Delta X^{k+1, m}(t)|^p]$ is finite from (\ref{est1}), we deduce
\begin{align*}
\lim_{N\rightarrow+\infty}\sup_{0\leq t\leq T}\bar{\mathbb{E}}\bigg[\bigg|\frac{\rho_2(|\Delta X^{k+1, m}(t)|^2)}{F_\varepsilon(\Delta X^{k+1, m, N}(t))}-\frac{\rho_2(|\Delta X^{k+1, m}(t)|^2)}{F_\varepsilon(\Delta X^{k+1, m}(t))}\bigg|\bigg]=0.
\end{align*}
Due to (\ref{est7}) again, the left-hand side of (\ref{new89}) converges to $\bar{\mathbb{E}}[F_\varepsilon(\Delta X^{k+1, m}(t))]$, as $N\rightarrow +\infty$. Then, by the Fatou-Lebesgue theorem, we have
\begin{align*}
\bar{\mathbb{E}}[F_\varepsilon(\Delta X^{k+1, m}(t))]
&\leq C\limsup_{N\rightarrow +\infty}\int_0^t\beta(s)\bigg(\rho_1(\bar{\mathbb{E}}[|\Delta X^{k, m}(s)|])+\bar{\mathbb{E}}\bigg[\frac{\rho_2(|\Delta X^{k+1, m}(s)|^2)}{F_\varepsilon(\Delta X^{k+1, m, N}(s))}\bigg]\bigg)ds\\
&\leq C\int_0^t\beta(s)\bigg(\rho_1(\bar{\mathbb{E}}[|\Delta X^{k, m}(s)|])+\limsup_{N\rightarrow +\infty}\bar{\mathbb{E}}\bigg[\frac{\rho_2(|\Delta X^{k+1, m}(s)|^2)}{F_\varepsilon(\Delta X^{k+1, m, N}(s))}\bigg]\bigg)ds\\
&= C\int_0^t\beta(s)\bigg(\rho_1(\bar{\mathbb{E}}[|\Delta X^{k, m}(s)|])+\bar{\mathbb{E}}\bigg[\frac{\rho_2(|\Delta X^{k+1, m}(s)|^2)}{F_\varepsilon(\Delta X^{k+1, m}(s))}\bigg]\bigg)ds.
\end{align*}
Letting $\varepsilon\rightarrow 0$, $F_\varepsilon(\Delta X^{k+1, m}(t))\downarrow |\Delta X^{k+1, m}(t)|$.
By Remark \ref{remconna}, for $t\in[0,T]$, $\lambda$-a.e., $\Delta X^{k+1, m}(t)$ belongs to $L^p_G(\Omega_t)$. One the other hand, for each $\varepsilon>0$, $F_
\varepsilon(x)$ is Lipschitz in $x$, then $F_\varepsilon(\Delta X^{k+1, m}(t))$ is also an element in $L^p_G(\Omega_t)$. By Theorem \ref{CTD}, (H\ref{H1'}') and Lemma \ref{le 2.2.2}, we obtain 
$\bar{\mathbb{E}}[F_\varepsilon(\Delta X^{k+1, m}(t))]\downarrow \bar{\mathbb{E}}[|\Delta X^{k+1, m}(t)|]$ and the following inequality:
\begin{align}\label{abcde}
\bar{\mathbb{E}}[|\Delta X^{k+1, m}(t)|]
&\leq C\int_0^t\beta(s)\bigg(\rho_1(\bar{\mathbb{E}}[|\Delta X^{k, m}(s)|])+\bar{\mathbb{E}}\bigg[\frac{\rho_2(|\Delta X^{k+1, m}(s)|^2)}{F_\varepsilon(\Delta X^{k+1, m}(s))}\bigg]\bigg)ds\notag\\
&\leq C\int_0^t\beta(s)\bigg(\rho_1(\bar{\mathbb{E}}[|\Delta X^{k, m}(s)|])+\bar{\mathbb{E}}\bigg[\frac{\rho_2(|\Delta X^{k+1, m}(s)|^2)}{|\Delta X^{k+1, m}(s)|}\bigg]\bigg)ds\\
&= C\int_0^t\beta(s)(\rho_1(\bar{\mathbb{E}}[|\Delta X^{k, m}(s)|])+\rho_3(\bar{\mathbb{E}}[|\Delta X^{k+1, m}(s)|]))ds.\notag
\end{align}
Borrowing the notation in the proof of Theorem \ref{th 3.1}, we rewrite (\ref{abcde}) into a simpler form:
$$u_{k+1,m}(t)\leq C\int_0^t\beta(s)(\rho_1(u_{k,m}(s))+\rho_3(u_{k+1,m}(s)))ds.$$
Taking the supremum of the left-hand side over all $m\in\mathbb{N}$, we have 
$$0\leq v_{k+1}(t)\leq C \int_0^t\beta(s)(\rho_1(v_k(s))+\rho_3(v_{k+1}(s)))ds,$$
and it follows that
$$0\leq \alpha(t)\leq C\int_0^t\beta(s)(\rho_1(\alpha(s))+\rho_3(\alpha(s)))ds.
$$
By (H\ref{H1'}') and Lemma \ref{le 2.2.1}, we deduce $$\alpha(t)=0,\ 0\leq t\leq T,$$
which implies that $\{X^m(\cdot)\}_{m\in\mathbb{N}}$ is a Cauchy sequence under the norm $\sup\limits_{0\leq t\leq T}\bar{\mathbb{E}}[|\cdot|]$. Hence, one can find a process $X$ in $M_G^1([0,T]; \mathbb{R}^n)$ such that
$$
\sup_{0\leq t\leq T}\bar{\mathbb{E}}[|X^m(t)-X(t)|]\rightarrow 0,\ {\rm{as}}\ N\rightarrow +\infty.
$$
By a classical argument, there exists a process $X(\cdot)\in M_G^1([0,T]; \mathbb{R}^n)$ and a subsequence
$\{X^{m_l}(\cdot)\}_{l\in\mathbb{N}}\subset\{X^m(\cdot)\}_{m\in\mathbb{N}}$ such that for each $t\in[0, T]$, 
$$X^{m_l}(t)\rightarrow X(t),\ {\rm{as}}\ l\rightarrow +\infty,\ {\rm{q.s.}}.$$
By the a priori estimate (\ref{est1})
and Lemma \ref{le 2.2.12}, we know
\begin{align*}
\sup_{0\leq t\leq T}\bar{\mathbb{E}}[|X^{m}(t)-X(t)|^p]&=\sup_{0\leq t\leq T}\bar{\mathbb{E}}[\liminf_{l\rightarrow +\infty}|X^m(t)-X^{m_l}(t)|^p]\\
&\leq \liminf_{l\rightarrow +\infty}\sup_{0\leq t\leq T}\bar{\mathbb{E}}[|X^{m}(t)-X^{m_l}(t)|^p]\leq CM<+\infty.
\end{align*}
Fixing a $\delta\in(0,+\infty)$, we calculate
\begin{align}\label{est13}
\limsup_{m\rightarrow+\infty}\sup_{0\leq t\leq T}\bar{\mathbb{E}}[|X^m(t)-X(t)|^2]
&\leq\delta^2+\limsup_{m\rightarrow+\infty}\sup_{0\leq t\leq T}\bar{\mathbb{E}}[|X^m(t)-X(t)|^2{\bf{1}}_{\{|X^m(t)-X(t)|>\delta\}}]\\
\notag&\leq\delta^2+\limsup_{m\rightarrow+\infty}\sup_{0\leq t\leq T}((\bar{\mathbb{E}}[|X^m(t)-X(t)|^p])^\frac{2}{p}(\bar{\mathbb{E}}[|{\bf{1}}_{\{|X^m(t)-X(t)|> \delta\}}|^\frac{p}{p-2}])^\frac{p-2}{p})\\
\notag&\leq\delta^2+M^\frac{2}{p}\limsup_{m\rightarrow+\infty}\sup_{0\leq t\leq T}(\bar{\mathbb{E}}[{\bf{1}}_{\{|X^m(t)-X(t)|> \delta\}}])^\frac{p-2}{p}\\
\notag&=\delta^2.
\end{align}
Because 
$$\lim_{m\rightarrow+\infty}\sup_{0\leq t\leq T}\bar{\mathbb{E}}[|X^m(t)-X(t)|]=0,$$
the last equality in (\ref{est13}) can be easily deduced by Lemma \ref{markov}. Letting $\delta\rightarrow 0$, we obtain
$$\lim_{m\rightarrow+\infty}\sup_{0\leq t\leq T}\bar{\mathbb{E}}[|X^m(t)-X(t)|^2]=0.$$
On the other hand, fixing a $\delta\in(0,+\infty)$, we have the following inequality in a similar way to (\ref{est18}):
\begin{align*}
\limsup_{m\rightarrow+\infty}\int_0^T\beta^2(t)\bar{\mathbb{E}}[|\rho_1(|X^{m}(t)&-X(t)|)|^2]dt\\
&\leq C(|\rho_1(\delta^2)|^2+K_\delta\lim_{m\rightarrow +\infty}\sup_{0\leq t\leq T}\bar{\mathbb{E}}[|X^{m}(t)-X(t)|^2])
\\&= C|\rho_1(\delta^2)|^2.
\end{align*}
As $\delta$
can arbitrary small, by (H\ref{H1'}'), we deduce
$$\lim_{m\rightarrow\infty}\bar{\mathbb{E}}\bigg[\sup_{0\leq t\leq
T}\bigg|\int_0^t(b(s,X^{m}(s))-b(s,X(s)))ds\bigg|^2\bigg]=0$$ 
and
$$\lim_{m\rightarrow\infty}\bar{\mathbb{E}}\bigg[\sup_{0\leq t\leq
T}\bigg|\int_0^t(h_{ij}(s,X^{m}(s))-h_{ij}(s,X(s)))d\langle
B^i,B^j\rangle_s\bigg|^2\bigg]=0,\ i, j = 1,\ldots,d.$$
Moreover, by the BDG type inequalities and Lemma \ref{le 2.2.2}, we can also deduce
\begin{align*}
\limsup_{m\rightarrow+\infty}\bar{\mathbb{E}}\bigg[\sup_{0\leq t\leq
T}\bigg|\int_0^t(g_j(s,X^{m}(s))-g_j(s,X(s)))dB^j_s\bigg|^2\bigg]
=0,\ j=1,\ldots, d.
\end{align*}
From all above, we conclude that $X\in M^2_G([0,T];\mathbb{R}^n)$ is a solution to the GSDE (\ref{eq 12}).\\[6pt]
Now we turn to the proof of uniqueness. Suppose $X_1$, $X_2\in M^2_G([0,T];\mathbb{R}^n)$ are two solutions that satisfy the GSDE (\ref {eq 12}), borrowing the notation in the proof of existence, we define for each $N\in\mathbb{N}$,
\begin{align*}
(X^1)^N(t)=x&+\int_0^tb^N(s,X^1(s))ds\\
&+\sum^d_{i,j=1}\int_0^th^N_{ij}(s,X^1(s))d\langle
B^i, B^j\rangle_s
+\sum^d_{j=1}\int_0^tg^N_j(s,X^1(s))d B^j_s,\ 0\leq t\leq T;
\end{align*}
\begin{align*}
(X^2)^N(t)=x&+\int_0^tb^N(s,X^2(s))ds\\
&+\sum^d_{i,j=1}\int_0^th^N_{ij}(s,X^2(s))d\langle
B^i, B^j\rangle_s
+\sum^d_{j=1}\int_0^tg^N_j(s,X^2(s))d B^j_s,\ 0\leq t\leq T.
\end{align*}
Following a similar procedure in the proof of existence, we know that
$\{(X^1)^N\}_{N\in\mathbb{N}}$ and $\{(X^2)^N\}_{N\in\mathbb{N}}$ converge to $X^1$ and $X^2$, respectively in $M^1_G([0,T]; \mathbb{R}^n)$ and we have
\begin{align*}
\bar{\mathbb{E}}[F_\varepsilon((X^1)^{N}(t)&-(X^2)^{N}(t))]\\
&\leq C\int_0^t\beta(s)\bigg(\rho_1(\bar{\mathbb{E}}[|X^1(s)-X^2(s)|])+
\bar{\mathbb{E}}\bigg[\frac{\rho_2(|X^1(s)-X^2(s)|^2)}{F_\varepsilon((X^1)^{N}(s)-(X^2)^{N}(s))}\bigg]ds.
\end{align*}
Letting $N\rightarrow+\infty$ and $\varepsilon\rightarrow 0$, we deduce
\begin{align*}
\bar{\mathbb{E}}[|X^1(t)-X^2(t)|]
&\leq C\int_0^t\beta(s)(\rho_1(\bar{\mathbb{E}}[|X^1(s)-X^2(s)|])+\rho_3(\bar{\mathbb{E}}[|X^1(s)-X^2(s)|]))ds.
\end{align*}
Thus,
$$
\sup_{0\leq s\leq t}\bar{\mathbb{E}}[|X^1(s)-X^2(s)|]\leq C\int_0^t\beta(s)(\rho_1+\rho_3)(\sup_{0\leq u\leq s}\bar{\mathbb{E}}[|X^1(u)-X^2(u)|])ds.
$$
Finally, Lemma \ref{le 2.2.1} gives the uniqueness result.\hfill{}$\square$
\begin{remark} Fang and Zhang \cite{FZ} have proved a pathwise uniqueness result for the classical SDEs by a stopping time technique, where $\rho$ is not necessary to be concave. Although we do have a similar stopping time technique, Lemma 3.3 in Fang and Zhang \cite{FZ} could not hold true in the $G$-framework, because for an $M^2_G((0, T);\mathbb{R}^n)$ process $\xi$, it is difficult to verify whether $\Phi(\xi)$ (using the notation in that paper) satisfies Definition 4.4 in Li and Peng \cite{LP} or not, that means the $G$-stochastic integrals in the proof of that lemma, whose upper limit involves a stopping time, could not be well defined. Fang and Zhang \cite{FZ} have also derived an existence result by the well-known Yamada-Watanabe theorem, which says that the existence of weak solution and pathwise uniqueness imply the existence of strong solution. In the $G$-framework, the corresponding Yamada-Watanabe theorem are unfortunately not available.
\end{remark}
\section{Solvability of $G$-backward stochastic differential equations}
\noindent In this section, we prove the existence and uniqueness theorem for the following GBSDE:
\begin{equation}\label{eq 15}
Y_t=\mathbb{E}\bigg[\xi+\int_t^Tf(s,Y_s)ds+\sum^d_{i,j=1}\int_t^Th_{ij}(s,Y_s)d\langle
B^i,B^j\rangle_s\bigg|\Omega_t\bigg],\ 0\leq t\leq T;
\end{equation}
where
$\xi\in L_G^1(\Omega_T; \mathbb{R}^n)$ and $f$, $g_{ij}$ are given functions that satisfy for each $x\in \mathbb{R}^n$, $f(\cdot, x)$, $h_{ij}(\cdot, x)\in M^1_G([0,T];\mathbb{R}^n)$, $i$, $j=1, \ldots, d$.\\[6pt]
We assume moreover that, for each $t\in [0, T]$ and $y$, $y_1$, $y_2\in\mathbb{R}^n$:
\begin{slist}
\item\label{H1"}$|f(t,y_1)-f(t,y_2)|+|h(t,y_1)-h(t,y_2)|\leq|\beta(t)|\rho(|y_1-y_2|)$;
\item\label{H2"}$|f(t,y)|+|h(t,y)|\leq\beta_1(t)+\beta_2|y|$,
\end{slist}
where $\beta:[0,T]\rightarrow\mathbb{R}^+$ is Lebesgue integrable, $\beta_1\in M_G^1([0,T])$, $\beta_2\in\mathbb{R}^+$ and  $\rho:(0,+\infty)\rightarrow(0,+\infty)$ is a continuous
increasing and concave function that vanishes at $0+$ and satisfies (\ref{eq 11}).
\begin{theorem}\label{th 3.3}
Under the assumptions above, (\ref{eq 15}) admits a unique solution
$Y\in M_G^1([0,T];\mathbb{R}^n).$
\end{theorem}
\noindent\textbf{Proof:} Let $Y_1,Y_2\in M_G^1([0,T];\mathbb{R}^n)$ be two solutions
of (\ref{eq 15}), then
\begin{align*}
Y_t^1-Y_t^2&=\mathbb{E}\bigg[\xi+\int_t^Tf(s,Y^1_s)ds+\sum^d_{i,j=1}\int_t^Th_{ij}(s,Y^1_s)d\langle
B^i,B^j\rangle_s\bigg|\Omega_t\bigg]\\
&-\mathbb{E}\bigg[\xi+\int_t^Tf(s,Y^2_s)ds+\sum^d_{i,j=1}\int_t^Th_{ij}(s,Y^2_s)d\langle
B^i,B^j\rangle_s\bigg|\Omega_t\bigg].
\end{align*}
Due to the sub-additivity of $\mathbb{E}[\cdot|\Omega_t]$, we obtain
\begin{align*}
|Y_t^1-Y_t^2|
&\leq\mathbb{E}\bigg[\bigg|\int_t^T(f(s,Y^1_s)-f(s,Y^2_s))ds\bigg|\bigg|\Omega_t\bigg]\\
&+\sum^d_{i,j=1}\mathbb{E}\bigg[\bigg|\int_t^T(h_{ij}(s,Y^1_s)-h_{ij}(s,Y^2_s))d\langle B^i,B^j\rangle_s\bigg|\bigg|\Omega_t\bigg].
\end{align*}
Taking $G$-expectation on both sides and using the BDG type inequalities and Lemma \ref{le 2.2.2}, we have 
\begin{align*}
\mathbb{E}[|Y_t^1-Y_t^2|]
&\leq\mathbb{E}\bigg[\bigg|\int_t^T(f(s,Y^1_s)-f(s,Y^2_s))ds\bigg|\bigg]\\
&+\sum^d_{i,j=1}\mathbb{E}\bigg[\bigg|\int_t^T(h_{ij}(s,Y^1_s)-h_{ij}(s,Y^2_s))d\langle
B^i,B^j\rangle_s\bigg|\bigg]\\
&\leq C\int_t^T \rho(\mathbb{E}[|Y_s^1-Y_s^2|])ds.
\end{align*}
Set
$$u(t)=\mathbb{E}[|Y_t^1-Y_t^2|],$$ then
$$u(t)\leq K\int_t^T \rho(u(s))ds.$$
By Lemma \ref{le 2.2.1}, we deduce
$$u(t)=0,\ 0\leq t\leq T,$$
which yields the pathwise uniqueness.\\[6pt]
For the proof of existence, we define a sequence of processes $\{Y^m\}_{m\in\mathbb{N}}$ as follows:
$$Y^0(t)=0,\ 0\leq t\leq T,$$
and 
$$Y^{m+1}_t=\mathbb{E}\bigg[\xi+\int_t^Tf(s,Y^m_s)ds+\sum^{d}_{i,j=1}\int_t^Th_{ij}(s,Y^m_s)d\langle
B^i,B^j\rangle_s\bigg|\Omega_t\bigg],\ 0\leq t\leq T.$$ The rest of the proof
goes in a similar way to the proof of Theorem \ref{th 3.1}, so we omit
it.\hfill{}$\square$
\begin{remark}
We notice that the definition of the GBSDE above is not the typical one (cf. (3.1) in Hu et al. \cite{HJPS}), in which the generator $f$ involves $Z$, i.e., the integrand of the It\^o type $G$-stochastic integral with respect to $G$-Brownian motion.
Based on the great efforts of many authors, such as Xu and Zhang \cite{XZ}, Soner et al. \cite{STZ2} and Song \cite{SS, S}, Peng et al. \cite{PSZ} have given a complete theory for $G$-martingale representation. Subsequently, Hu et al. \cite{HJPS} have derived a complete existence and uniqueness theorem for nonlinear GBSDEs with a generator $f$ that is uniformly Lipschitz in both $y$ and $z$.\\[6pt]
\noindent An extensive study to GBSDEs is meaningful because there will be numerous possible applications of GBSDEs in finance, for example, pricing and robust utility maximization in a model with a non-dominated class of probability measures.
\end{remark}
\section{Appendix}
\noindent
In the appendix, we give the proof of Lemma \ref{complex} in three steps. First of all, we consider the simplest case when $\zeta$ is uniformly Lipschitz in $x$. Then, we prove this lemma for all $\zeta$ that is uniformly bounded. To generalize the result to the case that $\zeta$ is unbounded, we need to define a sequence of truncated functions 
$\{\zeta^N\}_{N\in\mathbb{N}}$ as (\ref{trunccoe}) and complete the proof with the help of Theorem \ref{lpg}.
Now, we begin with the following lemmas.
\begin{lemma}\label{cs}
For some $p\geq 1$, $\zeta$ is a function that satisfies 
$\zeta(\cdot, x)\in M^p_G([0,T];\mathbb{R}^n)$ for each $x\in \mathbb{R}^n$. We assume moreover that $\zeta(\cdot, x)$ satisfies the Lipschitz condition, i.e., for each $t\in[0, T]$ and each $x_1$, 
$x_2\in \mathbb{R}^n$, $|\zeta(t, x_1)-\zeta(t, x_2)|\leq C_L|x_1-x_2|$.
Then, for each $X\in M^p_G([0,T];\mathbb{R}^n)$, $\zeta(\cdot, X_\cdot)$ is an element in $M^p_G([0,T];\mathbb{R}^n)$.
\end{lemma}
\noindent \textbf{Proof:} Without loss of the generality, we only give the proof of the one dimensional case. Suppose that $X$ can be approximated by a sequence $\{X^N\}_{N\in\mathbb{N}}\subset M^{p, 0}_G([0, T])$ of the form below:
$$
X^N_t:=\sum^{N-1}_{k=0}\xi_k{\bf{1}}_{[t_k, t_{k+1})}(t), 
$$
where $\xi_k\in L^0_{ip}(\Omega_{t_k})$, 
then
$$\int^T_0\bar{\mathbb{E}}[|\zeta(t, X^N_t)-\zeta(t, X_t)|^p]dt\leq C_L\int^T_0\bar{\mathbb{E}}[|X^N_t-X_t|^p]dt\rightarrow 0,\ {\rm{as}}\ N\rightarrow +\infty.$$
To obtain the desired result, we only need to prove that for each $k\in\mathbb{N}$, $\zeta(\cdot, \xi_k){\bf{1}}_{[t_k, t_{k+1})}(\cdot)\in M^p_G([0, T])$. In order to simplify the notation, we make a new assertion that is equivalent to the one stated above: fixing a $T\geq 1$, $\eta$ is an element in $L^0_{ip}(\Omega_1)$, then $\zeta(\cdot, \eta){\bf{1}}_{[1, T)}(\cdot)\in M^p_G([0, T])$. In what follows, we prove this assertion.\\[6pt]
Since $\eta\in L^0_{ip}(\Omega_1)$, there exists an $M>0$, such that $\eta\in [-M, M]$, which is a compact set in $\mathbb{R}$. For each $n\in\mathbb{N}$, we can find an open cover $\{G_i\}_{i\in I}$ of $\mathbb{R}$, such that $\lambda(G_i)<\frac{1}{n}$, $i\in I$. By the partition of unity theorem, there exists a family of $\mathcal{C}^\infty_0(\mathbb{R})$ function $\{\phi^n_i\}_{i\in I}$ such that for each $i\in I$, supp$(\phi^n_i)\in G_i$, $0\leq \phi^n_i \leq 1$ and for each $x\in\mathbb{R}$, $\sum\limits_{i\in I} \phi^n_i(x)=1$. Moreover, there exists a finite number of $\phi^n_i$ such that 
for $x\in[-M, M]$, $\sum\limits^{N(n)}_{i=1}\phi^n_i(x)=1$. Choosing for each $i= 1, \ldots, N(n)$ a point $x^n_i$ such that $\phi^n_i(x^n_i)>0$, we set 
$$\zeta^n(t, x)=\sum\limits^{N(n)}_{i=1}\zeta(t, x^n_i)\phi^n_i(x).$$ 
Then, $$
|\zeta^n(t, \eta){\bf{1}}_{[1, T)}(t)-\zeta(t, \eta){\bf{1}}_{[1, T)}(t)|\leq \sum^{N(n)}_{i=1}|\zeta(t, \eta)-\zeta(t, x^n_i)|\phi^n_i(\eta) \leq \frac{C_L}{n},\ 1\leq t < T,$$
which implies that $\zeta^n(\cdot, \eta){\bf{1}}_{[1, T)}(\cdot)$ converges to $\zeta(\cdot, \eta){\bf{1}}_{[1, T)}(\cdot)$ under the $M^p_G([0, T])$ norm (\ref{mnorm}). If for all $n\in\mathbb{N}$, $\zeta^n(\cdot, \eta){\bf{1}}_{[1, T)}(\cdot)$ belongs to $M^p_G([0, T])$, then, by the completeness of $M^p_G([0, T])$, $\zeta(\cdot, \eta){\bf{1}}_{[1, T)}(\cdot)\in M^p_G([0, T])$. It suffices to prove that $\zeta(\cdot, x^n_i)\phi^n_i(\eta){\bf{1}}_{[1, T)}(\cdot)\in M^p_G([0, T])$, $i= 1, \ldots, N(n)$,
which is given by the following lemma.\hfill$\square$
\begin{lemma}
Fixing a $T\geq 1$, let $X$ be an element in $M^p_G([0, T])$ and $\eta$ is an element in $L^0_{ip}(\Omega_1)$, then $\eta X_\cdot{\bf{1}}_{[1, T)}(\cdot)\in M^p_G([0, T])$ .
\end{lemma}
\noindent \textbf{Proof:} Suppose $X$ can be approximated by a sequence $\{X^N\}_{N\in\mathbb{N}}\subset M^{p, 0}_G([0, T])$ of the form below:
$$
X^N_t:=\sum^{N-1}_{k=0}\xi_k{\bf{1}}_{[t_k, t_{k+1})}(t), 
$$
then $X_\cdot{\bf{1}}_{[1, T)}(\cdot)$ can be approximated by a sequence $\{\bar{X}^N\}_{N\in\mathbb{N}}\subset M^{p, 0}_G([0, T])$:
$$
\bar{X}^N_t:=\sum^{N-1}_{k=0}\xi_k{\bf{1}}_{[t_k\vee 1, t_{k+1}\vee 1)}(t), 
$$
where $\xi_k\in L^0_{ip}(\Omega_{t_k})$. We define a sequence $\{\tilde{X}^N\}_{n\in\mathbb{N}}$ by
$$
\tilde{X}^N_t:=\sum^{N-1}_{k=0}\alpha_k{\bf{1}}_{[t_k, t_{k+1})}(t), 
$$
where
$$
\alpha_k:=
\left\{\begin{array}{c@{,}l}
0 & \ {\rm{if}}\ t_{k+1}< 1;\\[3pt]
\eta\xi_k &\ {\rm{if}}\ t_{k+1}\geq 1.
\end{array}\right.
$$
Since $L^0_{ip}(\Omega_1)\subset L^0_{ip}(\Omega_{t\vee 1})$ and $L^0_{ip}(\Omega_{t\vee 1})$ is closed under multiplication, we deduce that  $\{\tilde{X}^N\}_{N\in\mathbb{N}}\subset M^{p, 0}_G([0, T])$. Moreover, 
$$
|\tilde{X}^N_t-\eta X_t{\bf{1}}_{[1, T)}(t)|\leq |\eta\bar{X}^N_t-\eta X_t{\bf{1}}_{[1, T)}(t)|\leq M |X^N_t-X_t|,\ 0\leq t < T,
$$
where $M$ is the bound of $\eta$. This implies that $\eta X_\cdot{\bf{1}}_{[1, T)}(\cdot)$ is the limit of $\tilde{X}^N$ under the $M^p_G([0, T])$ norm (\ref{mnorm}).\hfill$\square$\\[6pt] 
\noindent\textbf{Proof of Lemma \ref{complex}:} 
Let $J\in \mathcal{C}^\infty(\mathbb{R}^n)$ be a non-negative function satisfies supp$(J)\subset B(0,1)$ and
$$
\int_{\mathbb{R}^n}J(x)dx=1.
$$
For each $\lambda>0$, we set
$$J_\lambda(x)=\frac{1}{\lambda^n}J(\frac{x}{\lambda})$$
and
$$\zeta_\lambda(t, x)=\int_{\mathbb{R}^n} J_\lambda (x-y) \zeta(t,y) dy.$$
We assume that $\zeta$ is uniformly bounded, then $\zeta_\lambda$ is uniformly Lipschitz in $x$.
By Lemma \ref{cs}, we have $\zeta_\lambda(\cdot, X_\cdot)\in M^q_G([0,T];\mathbb{R}^n)$. To deduce the desired result, we only need to show that $\zeta(\cdot, X_\cdot)$ is the limit of $\zeta_\lambda(\cdot, X_\cdot)$ under the $M^q_G([0,T];\mathbb{R}^n)$ norm (\ref{mnorm}).\\[6pt]
Fixing a $\lambda>0$, we calculate
\begin{align*}
|\zeta_\lambda(t, x)-\zeta(t, x)|\leq\int_{\mathbb{R}^n}J_\lambda(y)|\zeta(t,x-y)-\zeta(t,x)|dy.
\end{align*}
Therefore, 
\begin{align*}
\int^T_0\bar{\mathbb{E}}[|\zeta_\lambda(t, X(t))-\zeta(t, X(t))|^q]dt&\leq\int^T_0\bar{\mathbb{E}}\bigg[\bigg|\int_{\mathbb{R}^n}J_\lambda(y)|\zeta(t,X(t)-y)-\zeta(t,X(t))|dy\bigg|^q\bigg]dt\\
&\leq\int^T_0|\beta(t)|^q\bigg(\int_{\mathbb{R}^n}J_\lambda(y)\gamma(|y|)dy\bigg)^qdt\\
&\leq|\gamma(\lambda)|^q\int^T_0|\beta(t)|^qdt\leq C|\gamma(\lambda)|^q\rightarrow 0,\ {\rm{as}}\  \lambda\rightarrow 0.
\end{align*}
For an unbounded function $\zeta$, we construct a sequence of processes $(\zeta^N)_{n\in\mathbb{N}}$ as (\ref{trunccoe}). Fixing an $N\in\mathbb{N}$, we have
\begin{align}\label{howareyou}
\int_0^T\bar{\mathbb{E}}[|\zeta^N(t,X(t))-\zeta(t,X(t))|^q]dt
&\leq\int_0^T\bar{\mathbb{E}}[|\zeta(t,X(t))|^q{\bf{1}}_{\{|\zeta(t,X(t))|>N\}}]dt\\
&\leq\int_0^T\bar{\mathbb{E}}[(\beta_1(t)+\beta_2|X(t)|)^q{\bf{1}}_{\{\beta_1(t)+\beta_2|X(t)|>N\}}]dt\notag\\
&\leq C\bigg(\int_0^T\bar{\mathbb{E}}[|\beta_1(t)|^q{\bf{1}}_{\{|\beta_1(t)|>\frac{N}{2}\}}]dt
+\int_0^T\bar{\mathbb{E}}[(\beta_2|X(t)|)^q {\bf{1}}_{\{\beta_2|X(t)|>\frac{N}{2}\}}]dt\bigg)\notag\\
&\leq C\bigg(\int_0^T\bar{\mathbb{E}}[|\beta_1(t)|^q{\bf{1}}_{\{|\beta_1(t)|>\frac{N}{2}\}}]dt
+\beta_2^q\int_0^T\bar{\mathbb{E}}[|X(t)|^q {\bf{1}}_{\{|X(t)|>\frac{N}{2\beta_2}\}}]dt\bigg).\notag
\end{align}
Since $\beta_1$ and $X$ are $M^q_G([0,T])$ processes, by Remark \ref{remconna} and Theorem \ref{lpg} along with Lebesgue's dominated convergence theorem, the right-hand side of (\ref{howareyou}) converges to 0. This yields the desired result.\hfill{}$\square$\\[6pt]
\noindent\textbf{Acknowledgement} The authors express special thanks to Prof. Ying Hu, who provided both the initial inspiration for the work and useful suggestions.

\end{document}